\documentclass[12pt]{article}
\usepackage{amssymb,amsmath}
\usepackage{graphicx}
\usepackage{amsfonts}
\usepackage{float}
\usepackage{flafter}
\usepackage{color}

\labelwidth\leftmargini\advance\labelwidth-\labelsep

\def\R{\relax\ifmmode I\!\!R\else$I\!\!R$\fi}

\def\Z{\relax\ifmmode Z\!\!\!Z\else$Z\!\!\!Z$\fi}

\def\C{\relax\ifmmode C\!\!\!\!I\else$C\!\!\!\!I$\fi}

\def\K{\relax\ifmmode I\!\!K\else$I\!\!K$\fi}

\def\N{\relax\ifmmode I\!\!N\else$I\!\!N$\fi}

\newcounter{defcounter}[section]
{\vspace{0.1cm}\begin{sloppypar}\noindent\stepcounter{defcounter}{\bfseries
Definition
      \thesection.\thedefcounter}}%
{\end{sloppypar}\vspace{0.1cm}}

\newtheorem{theorem}{Theorem}[section]

\newtheorem{proposition}{Proposition}[section]

\newcommand{\proof}{{\bf Proof.} }

\newcommand{\qed}{\hfill $\square$}

\begin{document}
\thispagestyle{empty}
\begin{center}
{\Large {\bf Exceptions in the domain of generic absolute continuity of non-homogeneous self-similar measures}}
\end{center}
\begin{center}J\"org Neunh\"auserer\\
Technical University of Braunschweig \\
joerg.neunhaeuserer@web.de
\end{center}
\begin{center}
\begin{abstract}
Non-homogeneous self-similar measures are generically absolute continuous in the domain of parameters for which the similarity dimension is larger than one, see \cite{[SSS]}. Using certain algebraic curves we construct here exceptional singular non-homogeneous self-similar measures in this domain.  
\end{abstract}
\end{center}
\section{Introduction and main result}
For $\beta_{1},\beta_{2}\in(0,1)$ we consider the similarities on $\mathbb{R}$ given by
\[ T_{1}(x)=\beta_{1}x+\beta_{1}\mbox{ and } T_{2}(x)=\beta_{2}x-\beta_{2}.\]
It is well known that for probability $p\in(0,1)$ there is an unique self-similar Borel probability measure $\mu^{p}_{\beta_{1},\beta_{2}}$ on $\mathbb{R}$ with
\[ \mu^{p}_{\beta_{1},\beta_{2}}=pT_{1}(\mu^{p}_{\beta_{1},\beta_{2}})+(1-p)T_{-1}(\mu^{p}_{\beta_{1},\beta_{2}}),\]
see \cite{[HU]}. If $\beta_{1}=\beta_{2}=\beta\in(0,1)$ the measure $\mu^{p}_{\beta}=\mu^{p}_{\beta,\beta}$ is an infinit convolved Bernoulli measure. These measure were intensively studied since the pioneering work of Erdös \cite{[ER1],[ER2]} with milestone results of Solomyak \cite{[SO]}, Hochman \cite{[HO]}, Shmerkin \cite{[SCH]} and Varju \cite{[VA]}. In the case $\beta_{1}\not=\beta_{2}$ the measures $\mu^{p}_{\beta_{1},\beta_{2}}$ are now usually called non-homogeneous self-similar.\footnote{We were used to call these measure non-uniform self-similar but adopt the terminology from the decisive book \cite{[BA]}.} 
The similarity dimension of the function system $\{T_{1},T_{2}\}$, with respect to the probability vector $(p,(1-p))$, is defined by
\[ SD^{p}_{\beta_{1},\beta_{2}}=\frac{-p\log(p)-(1-p)\log(1-p)}{-p\log\beta_{1}-(1-p)\log\beta_{2}}.\]
The measures $\mu^{p}_{\beta_{1},\beta_{2}}$ are of pure type, they are either totally singular or absolutely continuous and even equivalent to the Lebesgue measure in the second case, see proposition 3.1 in \cite{[PS]}. Since the work of Hutchinson \cite{[HU]} we know 
that Hausdorff dimension of $\mu^{p}_{\beta_{1},\beta_{2}}$ is bounded by the similarity dimension
\[ \dim_{H}\mu^{p}_{\beta_{1},\beta_{2}}\le SD^{p}_{\beta_{1},\beta_{2}}.\]    
Hence the measures are singular if $SD^{p}_{\beta_{1},\beta_{2}}<1$. The main result on absolut continuity reads as follows:
\begin{theorem}
For all $p\in (0,1)$ and almost all $\beta_{1},\beta_{2}\in(0,1)$ with $SD^{p}_{\beta_{1},\beta_{2}}>1$ the measure $\mu^{p}_{\beta_{1},\beta_{2}}$ is absolutely continuous. 
\end{theorem}
We proved this result for $\beta_{1},\beta_{2}\in(0,0.649)$ in \cite{[NE1]} using transversality techniques which lead to the crude bound $0.649$. Ngai and Wang \cite{[NW]} proved a similar result using transversality. Saglietti, Shmerkin and Solomyak \cite{[SSS]} succeeded in removing the bound $0.649$ by an improvement of the classical tranversality method. For infinit convolved Bernoulli measure $\mu_{\beta}^{p}$ we have an even stronger result: The set of expectations to absolute continuity in the domain $\beta>p^p(1-p)^(1-p)$ has Hausdorff measure zero, see \cite{[SCH]}. The only known exceptions found by Erdös \cite{[ER1]} are algebraic numbers $\beta\in(0.5,1)$, which are reciprocals of Pisot numbers.\footnote{A Pisot numbers is an algebraic integer $\alpha>1$ with all its algebraic conjugates inside the unit circle.} In this case we have $\dim_{H}\mu_{\beta}^{p}<1$. In the non-homogeneous case we proved in \cite{[NE2]} the existence of  $\beta_{1}\not=\beta_{2}$ with $\beta_{1}\beta_{2}>1/4$ such that $\dim_{H}\mu^{0.5}_{\beta_{1},\beta_{2}}<1$, without construction explicit examples. Here we construct explicit algebraic curves with points $\beta_{1},\beta_{2}\in(0,1)$ such that $SD^{p}_{\beta_{1},\beta_{2}}>1$ but $\dim_{H}\mu^{p}_{\beta_{1},\beta_{2}}<1$ for some $p\in(0,1)$. These are exceptions in the domain of generic absolute continuity of non-homogeneous self-similar measures.\\
For $n\ge 3$ consider a finite sequence $s\in\{-1,1\}^{n}$. For $1\le k\le n$ let $\sharp_{k}(s)$ be the number of entries in $(s_{1},s_{2},\dots,s_{k})$ that are $1$ and $\tilde\sharp_{k}(s)=k-\sharp_{k}(s)$ the number of entries that are $-1$. For two sequences $s,t\in\{-1,1\}^{n}$ with $s\not=t$ and $\sharp_{n}(s)=\sharp_{n}(t)$ we define an algebraic curve $c_{s,t}$ by the algebraic equation  
\[ \sum_{k=1}^{\infty}s_{k}x^{\sharp_{k}(s)}y^{\tilde\sharp_{k}(s)}-t_{k}x^{\sharp_{k}(t)}y^{\tilde\sharp_{k}(t)}=0\]
in $\mathbb{R}^2$. Let 
\[\mathfrak{C}=\{c_{s,t}|s,t\in \{-1,1\}^{n},s\not=t, \sharp_{n}(s)=\sharp_{n}(t), n\ge 3\}\] be the set of all of such curves. With this notations we formulate our main result. 
\begin{theorem}
If $(\beta_{1},\beta_{2})\in(0,1)^2$ with $\beta_{1}+\beta_{2}>1$ is a point on a curve in $\mathfrak{C}$, then there is a $p\in(0,1)$ such that
 $SD^{p}_{\beta_{1},\beta_{2}}>1$, but the measure $\mu^{p}_{\beta_{1},\beta_{2}}$ is singular with $\dim_{H}\mu^{p}_{\beta_{1},\beta_{2}}<1$.
\end{theorem}
We prove this theorem in the next section. So far it is open if their are any curves in $\mathfrak{C}$ that have points in the open rectangle 
\[ R=\{(\beta_{1},\beta_{2})\in(0,1)^2|\beta_{1}+\beta_{2}>1\}.\]
For $n=3$ and $n=4$ some calculation show that there do not exists such curves. For $n=5$ we consider $s=(1,-1,-1,-1,1)$ and $t=(-1,1,1,-1,-1)$, which lead to the algebraic equation  
\[  2 x^2 y^3+x^2 y^2-x^2 y-x y^3-x y^2-2 x y+x+y =0.\]
The corresponding algebraic curve $c_{s,t,}$ has approbate points, see figure 1. For $n=5$ this is the only curve with the required properties found.
In the last section of the paper we will give suffizient conditions on the sequences $s,t\in\{-1,1\}^{n}$ such that $c_{s,t}\cap R\not=\emptyset$. These conditions especially apply to the example described here. Examples of curves which have the required properties with $n=6$ will be given there.     
\begin{figure}
\vspace{0pt}\hspace{40pt}\scalebox{0.35}{\includegraphics{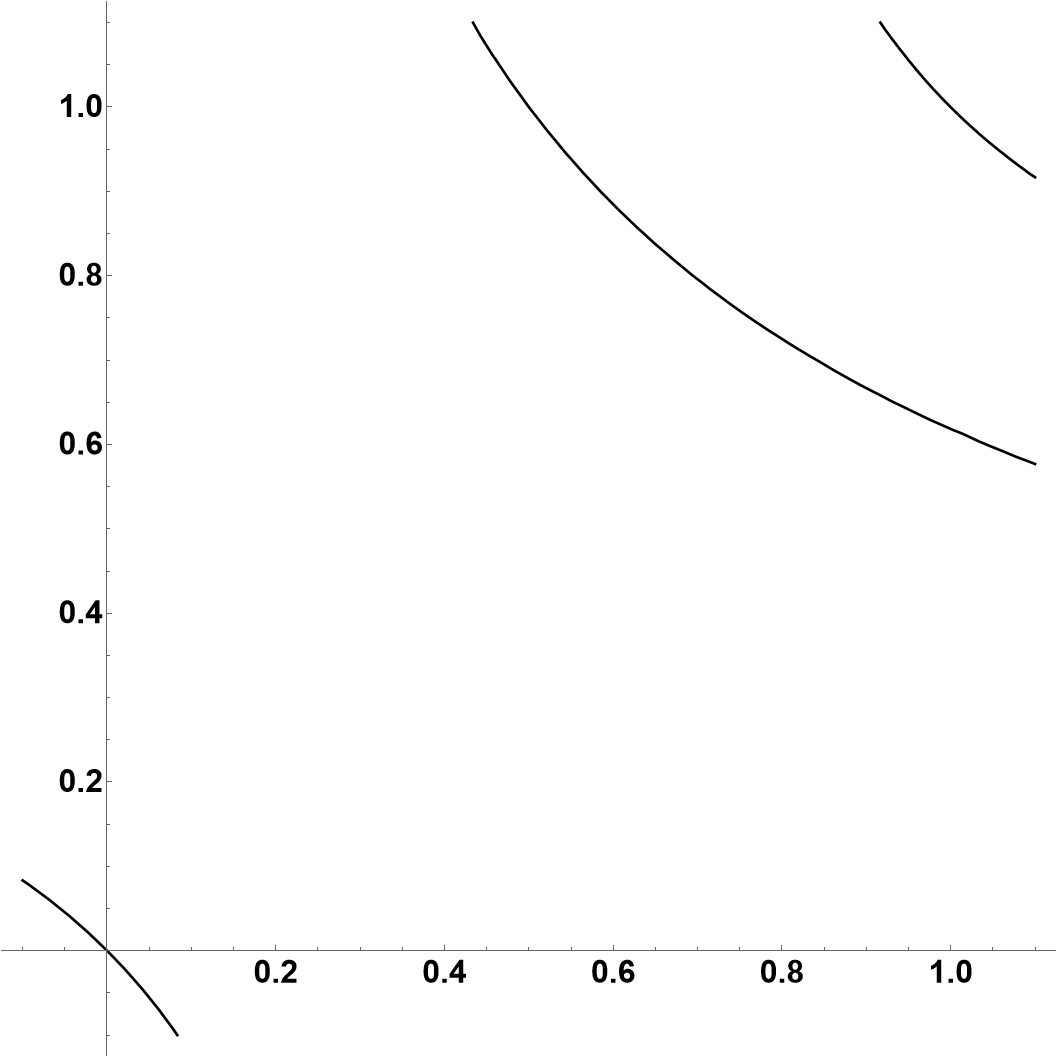}}
\caption{The curve $c_{s,t}$ for s=(1,-1,-1,-1,1) and $t=(-1,1,1,-1,-1)$ }
\end{figure}
\section{Proof of theorem 1.2}
We first estimate the Hausdorff dimensions of the measures $\mu^{p}_{\beta_{1},\beta_{2}}$ under the assumptions of theorem 1.2.
\begin{proposition}
If $(\beta_{1},\beta_{2})\in(0,1)^2$ is a point on a curve in $\mathfrak{C}$ than 
\[ \dim_{H}\mu^{p}_{\beta_{1},\beta_{2}}\le \widehat{SD}^{p}_{\beta_{1},\beta_{2}}< SD^{p}_{\beta_{1},\beta_{2}}\]     
for all $p\in(0,1)$ and $\widehat{SD}^{p}_{\beta_{1},\beta_{2}}$ is continuous in $p$.
\end{proposition}
 \proof For $r\in\{-1,1\}^{n}$ we set
 \[ T_{r}(x)=T_{r_{n}}\circ T_{r_{2}}\circ\dots \circ T_{r_{1}}(x)\]
 \[ = \beta_{1}^{\sharp_{n}(r)}\beta_{2}^{\tilde\sharp_{n}(r)}x+\sum_{k=1}^{\infty}r_{k}\beta_{1}^{\sharp_{k}(r)}\beta_{2}^{\tilde\sharp_{k}(r)}. \]
 Note that if $(\beta_{1},\beta_{2})\in(0,1)^2$ is a point on a curve $c_{s,t}\in\mathfrak{C}$ for $s,t\in\{-1,1\}^{n}$, than $T_{s}=T_{t}$.  
 By \cite{[HU]} there is a unique Borel probability measure $\mu$ on $\mathbb{R}$ fulfilling
 \[ \mu=\sum_{r\in \{-1,1\}^{n}}p^{\sharp_{n}(r)}(1-p)^{\tilde\sharp_{n}(r)}T_{r}(\mu)\]
 \[  = \sum_{r\in \{-1,1\}^{n}\backslash\{s,t\}}p^{\sharp_{n}(r)}(1-p)^{\tilde\sharp_{n}(r)}T_{r}(\mu)+ (p^{\sharp_{n}(s)}(1-p)^{\tilde\sharp_{n}(s)}+ p^{\sharp_{n}(t)}(1-p)^{\tilde\sharp_{n}(t)})T_{s}(\mu)\]
 and we obviously have $\mu=\mu^{p}_{\beta_{1},\beta_{2}}$. The similarity dimension of the function system $\{T_{r}|r\in\{-1,1\}^{n}\backslash\{s,t\}\}\cup\{T_{s}\}$ with respect to the probability vector that is summing up the probability of $s$ and $t$, 
 \[ p_{s,t}=(p^{\sharp_{n}(s)}(1-p)^{\tilde\sharp_{n}(s)}+ p^{\sharp_{n}(t)}(1-p)^{\tilde\sharp_{n}(t)})\] 
 is given by
 \[  \widehat{SD}^{p}_{\beta_{1},\beta_{2}}=\frac{-\sum_{r\in\{-1,1\}^{n}\backslash\{s,t\}}p^{\sharp_{n}(r)}(1-p)^{\tilde\sharp_{n}(r)}
 \log(p^{\sharp_{n}(r)}(1-p)^{\tilde\sharp_{n}(r)})-p_{s,t}\log(p_{s,t})}{-\sum_{r\in\{-1,1\}^{n}\backslash\{s,t\}}p^{\sharp_{n}(r)}(1-p)^{\tilde\sharp_{n}(r)}
 \log(\beta_{1}^{\sharp_{n}(r)}\beta_{2}^{\tilde\sharp_{n}(r)})-p_{s,t}\log(\beta_{1}^{\sharp_{n}(s)}\beta_{2}^{\tilde\sharp_{n}(s)})}\]  
 
 \[ =\frac{-\sum_{r\in\{-1,1\}^{n}\backslash\{s,t\}}p^{\sharp_{n}(r)}(1-p)^{\tilde\sharp_{n}(r)}
 \log(p^{\sharp_{n}(r)}(1-p)^{\tilde\sharp_{n}(r)})-p_{s,t}\log(p_{s,t})}{-p\log\beta_{1}-(1-p)\log\beta_{2}}\]
 \[<\frac{-p\log(p)-(1-p)\log(1-p)}{-p\log\beta_{1}-(1-p)\log\beta_{2}}=SD^{p}_{\beta_{1},\beta_{2}}.\]
The inequality here is due to
\[ -p_{s,t}\log(p_{s,t})< -(p^{\sharp_{n}(s)}(1-p)^{\tilde\sharp_{n}(s)}\log((p^{\sharp_{n}(s)}(1-p)^{\tilde\sharp_{n}(s)})-
p^{\sharp_{n}(t)}(1-p)^{\tilde\sharp_{n}(t)}\log(p^{\sharp_{n}(t)}(1-p)^{\tilde\sharp_{n}(t)}).\] 
The proposition now follows from \cite{[HU]}. 
 \qed~\\~\\
To complete the proof of theorem 1.2 it remains to show:
\begin{proposition}
If $(\beta_{1},\beta_{2})\in(0,1)^2$ with $\beta_{1}+\beta_{2}>1$ is a point on a curve in $\mathfrak{C}$ than there is an $p\in(0,1)$ such that $\widehat{SD}^{p}_{\beta_{1},\beta_{2}}<1$ and $SD^{p}_{\beta_{1},\beta_{2}}>1$. Here $\widehat{SD}^{p}_{\beta_{1},\beta_{2}}<1$ is given by the last proposition. 
\end{proposition} 
\proof Since $\beta_{1}+\beta_{2}>1$ there is an $d>1$ such that $\beta_{1}^{d}+\beta_{2}^{d}=1$. Let $p_{M}=\beta_{1}^{d}\in(0,1)$. A simple calculation shows
\[ SD^{p_{M}}_{\beta_{1},\beta_{2}}=d>1.\]
Since 
\[ \lim_{p\to 0}SD^{p}_{\beta_{1},\beta_{2}}=\lim_{p\to 1}SD^{p}_{\beta_{1},\beta_{2}}=0\] 
there is a $p_{1}\in(0,1)$ such that $SD^{p_{1}}_{\beta_{1},\beta_{2}}=1$. Since $\widehat{SD}^{p}_{\beta_{1},\beta_{2}}< SD^{p}_{\beta_{1},\beta_{2}}$ for all $p\in(0,1)$ by proposition 2.1 and $SD^{p}_{\beta_{1},\beta_{2}}$ and $\widehat{SD}^{p}_{\beta_{1},\beta_{2}}$ are continuous in $p$
there is an $p$ near $p_{1}$ such that $\widehat{SD}^{p}_{\beta_{1},\beta_{2}}<1$ and $SD^{p}_{\beta_{1},\beta_{2}}>1$.    
\qed
\section{Algebraic curves given the exceptions}
To apply theorem 2.1 we need sequences $s,t\in\{-1,1\}^{n}$ such that the algebraic curves $c_{s,t}\in\mathfrak{C}$, defined in the first section, have points in the open rectangle  
\[ R=\{(x,y)\in(0,1)^2|x+y>1\}. \]
We find a sufficient condition on the sequences that guarantees the existences of such points.
\begin{proposition}
For $n\ge 3$ let $s,t\in \{-1,1\}^{n}$ with $s\not=t$ and $\sharp_{n}(s)=\sharp_{n}(t)$. Assume that $s$ begins with $(1,-1)$ and $t$ begins with $(-1,1)$. If  
\[ \sum_{k=1}^{\infty}s_{k}\sharp_{k}(s)>\sum_{k=1}^{\infty}t_{k}\sharp_{k}(t),\]
we have $c_{s,t}\cap R\not=\emptyset$.
\end{proposition}  
\proof
Let
\[ f(x,y)=\sum_{k=1}^{\infty}s_{k}x^{\sharp_{k}(s)}y^{\tilde\sharp_{k}(s)}-t_{k}x^{\sharp_{k}(t)}y^{\tilde\sharp_{k}(t)}\]
and
\[ f(x)=f(x,1)=\sum_{k=1}^{\infty}s_{k}x^{\sharp_{k}(s)}-t_{k}x^{\sharp_{k}(t)}.\]
By the first assumption we have $f(0)=1$ and since $\sharp_{n}(s)=\sharp_{n}(t)$ we have $f(1)=0$. By the second assumption we obtain
\[ f^{\prime}(1)=\sum_{k=1}^{\infty}s_{k}\sharp_{k}(s)-\sum_{k=1}^{\infty}t_{k}\sharp_{k}(t)>0,\]
which implies $f(x)<0$ for some $x\in(0,1)$. It follows that there is an $x_{0}\in(0,1)$ and $\delta>0$ such that $f(x_{0},1)=f(x_{0})=0$ and $f(x)>0$ for $x\in(x_{0}-\delta,x_{0})$ and $f(x)<0$ for $x\in(x_{0},x_{0}+\delta)$. Hence by continuity for all suffizient small $\epsilon>0$ there exists a 
$x_{-}\in(x_{0}-\delta,x_{0})$ with $f(x_{-},1-\epsilon)<0$ and $x_{+}\in(x_{0},x_{0}+\delta)$ with $f(x_{+},1-\epsilon)<0$. Again we continuity there is an $x_{\epsilon}\in(x_{-},x_{+})$ with $f(x_{\epsilon},1-\epsilon)=$. Moreover if $\epsilon$ is suffizient small $(x_{\epsilon},1-\epsilon)\in R$.
\qed~\\\\
This proposition applies to the sequences $s,t$ with $n=5$ given in the first section of the paper. For $n=6$ we have the following examples:
\[  s=(1,-1,-1,1,1,1)\quad t=(-1,1,1,1,1,-1),\]
\[  s=(1,-1,1,-1,1,1)\quad t=(-1,1,1,1,1,-1),\]
\[  s=(1,-1,-1,-1,1,1)\quad t=(-1,1,-1,1,1,-1),\]
\[  s=(1,-1-,1,-1,1,1)\quad t=(-1,1,1,-1,1,-1),\]
given the algebraic equations
\[ 2 x^4 y^2 - x^4 y + x^3 y^2 - x^3 y + x^2 y^2 - x^2 y - x y^2 -2 x y + x + y=0,\]
\[ 2 x^4 y^2 - x^4 y + x^3 y^2 - x^3 y - x^2 y^2 - 2 x y + x + y=0,\]
\[ 2x^3 y^3 - x^3 y^2 + x^2 y^3 - x^2 y^2 - x y^3 - 2 x y + x + y=0, \]
\[ 2 x^3 y^3 - x^3 y^2 + x^2 y^3 + x^2 y^2 - x^2 y - x y^3 - x y^2 -2 x y + x + y=0.\] 
We do not know if the condition given in the last proposition is (up to symmetries) necessary to guarantee the existentes of appropriate curves. We leave this question to the reader.      

\end{document}